\documentclass{amsproc}

\usepackage{url}

\usepackage{cite}

    \usepackage{amsthm,amsmath,graphicx,color}

\usepackage{newunicodechar}

\theoremstyle{definition}

\theoremstyle{remark}

\numberwithin{equation}{section}

    \newcommand{\parrot}{\raisebox{-0.6ex}{\protect\includegraphics[height=3ex]{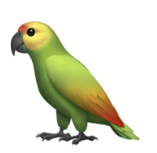}}}
    
\begin{document}

\title{A Non-Expert's Introduction to Data Ethics for Mathematicians}

\author{Mason A. Porter}
\address{Department of Mathematics, University of California, Los Angeles, Los Angeles, California 90095, USA; Department of Sociology, University of California, Los Angeles, Los Angeles, California 90095, USA; Santa Fe Institute, Santa Fe, New Mexico, 87501, USA}
\email{mason@math.ucla.edu}
\thanks{In parts of this chapter,
 I drew on materials from a class that Matt Salganik has taught on computational social science at Princeton University. Johan Ugander suggested several excellent resources and research studies. 
  I thank Jacob Foster, Mark Handcock, Andrei Kleshin, Peter Mucha, Dennis M\"uller, Deanna Needell, Matt Salganik, Massimo Stella, Alexandria Volkening, Chris Wiggins, Richard Yeh, and an anonymous referee for helpful comments.
  Finally, I thank my fellow Short Course organizers (Heather Zinn Brooks, Michelle Feng, and Alexandria Volkening) for their continuing support on this project and elsewhere, their collaboration, and their friendship.}

\date{January 18, 2022 and, in revised form, December 31, 2023.}

\dedicatory{This chapter is dedicated to my current and former Ph.D. students. They mean the world to me.}




\begin{abstract}
I give a short introduction to data ethics. I begin with some background information and societal context for data ethics. I then discuss data ethics in mathematical-science education and indicate some available course material. I briefly highlight a few efforts --- at my home institution and elsewhere --- on data ethics, society, and social good. I then discuss open data in research, research replicability and some other ethical issues in research, and the tension between privacy and open data and code, and a few controversial studies and reactions to studies. I then discuss ethical principles, institutional review boards, and a few other considerations in the scientific use of human data. I then briefly survey a variety of research and lay articles that are relevant to data ethics and data privacy. I conclude with a brief summary and some closing remarks.

My focal audience is mathematicians, but I hope that this chapter will also be useful to others. I am not an expert about data ethics, and this chapter provides only a starting point on this wide-ranging topic. I encourage you to examine the resources that I discuss and to reflect carefully on data ethics, its role in mathematics education, and the societal implications of data and data analysis. As data and technology continue to evolve, I hope that such careful reflection will continue throughout your life.
\end{abstract}

\maketitle



\begin{quotation}
\emph{Don't say that he's hypocritical \\
Say rather that he's apolitical \\
``Once the rockets are up, who cares where they come down? \\
That's not my department!'' says Wernher von Braun.} \\ \\
(Tom Lehrer, \emph{Wernher von Braun}, 1975)
\end{quotation}


\section{Introduction} 
\label{intro}

The use of digital data to examine and help understand human behavior is both powerful and dangerous \cite{editorial2021}. Every day, it seems like there is a new nightmare with problematic uses of data and algorithms. The use of predictive policing to identify criminal activity can exacerbate existing racial and ethnic inequities \cite{proactive2018}, algorithmic social-credit ratings of individuals have frightening dystopian uses \cite{social-credit}, and other manifestations of ``algocracy'' (i.e., algorithmic government) have many human ramifications \cite{engin2019,osoba2019,hassel2021,straub2022}. 

The ever-increasing use of and reliance on ``big data'' and the accelerating application of tools such as machine learning and artificial intelligence (AI) to make societally consequential decisions can cause very significant harm and exacerbate societal inequities \cite{cathy2016,royal2018}. Problems like algorithmic bias and the collection, measurement, and use of enormous amounts of data about humans and their behavior have significant societal consequences \cite{wagner2021,lazer2021,sadowski2021,birhane2021,wiggins2023}. New technologies like ``deep fakes'' \cite{ipam-fakery}, synthetic media in which a person in an existing image or video is replaced with somebody else's likeness,\footnote{A light-hearted example is Gollum's precious cover of the song ``Nothing Compares to U'' \cite{precious}.} have terrifying potential to cause harm. There are significant risks and potential nefarious uses of technology like generative AI and large language models (LLMs) \cite{parrots2021,mahdi2023,ferrara2023}. Modern data analysis also has been accompanied by the reincarnation of pseudosciences, such as physiognomy, which can now be performed at a large scale \cite{stark2023}. Along with such dangers come enormous potential benefits, and there are nascent efforts to create documents like an ``AI Bill of Rights" to help guide the design and deployment of automated systems in a way that protects people~\cite{OSTP2022}.

The situation is already scary. In 2016, Microsoft released a Twitter bot that ``learned'' from its interactions with other Twitter accounts; in less than a day, it was regularly producing racist tweets \cite{hatebot2016}. In December 2020, Stanford University's use of an algorithm to determine which people would receive the first batches of a COVID-19 vaccine resulted in an inequitable vaccination roll-out that prioritized high-ranking doctors over frontline human-facing medical personnel \cite{covid2020}. The output of algorithms have also led to decisions to cut off medical care \cite{ross2023}. In 2021, the outgoing Editor-in-Chief of a scientific journal publicly released refereeing data, including the numbers of decisions to accept or reject papers by each of the journal's referees during his tenure \cite{referee2021}. In 2023, the LLM ChatGPT invented\footnote{It seems that LLMs and generative AI tools often simply make things up.} a sexual-harassment scandal and named a real professor in it \cite{verma2023} and Microsoft limited its Bing AI chatbot after its unsettling conversations \cite{leswing2023}. In today's world, it often seems necessary to laugh to keep from crying \cite{smbc2021}.

The mathematical, statistical, and computational sciences are interconnected with human communities and society at large \cite{cambridge,ethicalmath,muller2023b,moses2021}. Neither our education nor our research occur in a vacuum, and many of our algorithms and other methods are now applied to social systems \cite{wagner2021,chiodo2023}. However, it is not traditional in mathematical-science education to discuss data ethics and other ethical considerations \cite{rycroft2022,muller2022}. What we choose to teach (and choose not to teach) impacts what our mentees do with their education. Consequently, it is crucial for mathematical scientists (and others) to think carefully about ethics and engage with it throughout their careers~\cite{chiodo2024}. See \cite{loukides2018} for a short book on ethics and data science, \cite{ethics2016} for a theme journal issue on data ethics, \cite{veliz2021} for a wide-ranging handbook of digital ethics, and \cite{hippo2021} for a recent discussion of the ethics landscape in mathematics and advocacy of a ``Hippocratic Oath'' for mathematics.

The present chapter, which I intend as an introductory resource about data ethics, is a companion to my oral presentation at the American Mathematical Society Short Course on Mathematical and Computational Methods for Complex Social Systems in January 2021~\cite{short2021}. My slides and video presentation are available at \cite{mason-slides}. I discuss many of the same key points in my presentation as in the present chapter, although my emphases often differ. I have also learned about more data-ethics resources in the last three years, and certain ethical issues and relevant technologies have gained increased prominence in that time. I encourage you to look at my slides and watch my presentation, and I especially encourage you to look at the resources that I discuss in those slides and in the present chapter. I am not an expert on data ethics --- and, to be frank,  writing this article has been accompanied by my most serious case of ``imposter syndrome'' in a long time --- and it is important that you look at what actual experts have to say. I will attempt to give some helpful thoughts and pointers to begin a journey in data ethics. It is important to continue to reflect on data ethics and the societal impact of data and data analysis throughout your life.

In this chapter, I cover a diverse variety of topics. I start with data ethics and education, which is important for almost the entire mathematical-science community. In Section \ref{education}, I discuss data ethics in the mathematics community and mathematical-science education. I highlight the importance of curricular data ethics in Section \ref{curriculum}, and I point to existing educational material in Section \ref{courses}. In Section \ref{ucla}, I describe some of our efforts at University of California, Los Angeles (UCLA) on data ethics and society. I also briefly mention a few efforts by others and point out an important caveat in efforts in ``data for social good". In Section \ref{open}, I discuss various ethical issues and tensions in research. Some of these topics, such as research replication and appropriately acknowledging others, are directly relevant to everybody in the mathematical-science community. Other topics are most directly relevant to people who use data in their research, although it is still important for others in the mathematical sciences to have some familiarity with them. I discuss research replication and open data and code in Section \ref{rep}, acknowledgements and licensing material for open use in Section \ref{license}, privacy concerns and its tension with open data and code in Section \ref{concern}, and controversies in studies and in reactions to studies in Section \ref{controversy}. In Section \ref{open-summary}, I summarize the key ideas of Section \ref{open}. In Section \ref{principle}, I discuss some ethical principles and other considerations in the use of human data in scientific research. In Section \ref{research}, I highlight a few research articles about issues in data ethics. I conclude in Section \ref{summary}.


\section{Data Ethics, the Mathematics Community, and Mathematical-Science Education} 
\label{education}

In this section, I advocate for data ethics in mathematical-science education and point to some existing course materials. Many students who obtain undergraduate or graduate degrees in the mathematical sciences will become data scientists or otherwise will work extensively with data. Therefore, data ethics needs to be a central part of their professional lives, and it is our duty as mathematics educators to mentor them to face these ethical challenges with careful reflection, wisdom, and humility.


\subsection{The Importance of Data Ethics in the Mathematics Community and Mathematical-Science Curricula}\label{curriculum}

The mathematical, statistical, and computational sciences are intertwined with human society \cite{cambridge,ethicalmath,muller2023b}. This notion is not new. Throughout his life, educator and civil-rights activist Bob Moses championed mathematical literacy as a civil right for all public-school children, with particular advocacy for those who are most vulnerable~\cite{moynihan2023}. In a 2021 article~\cite{moses2021} that was published a month after his death, Moses wrote
\medskip
\begin{quotation}
\emph{In the 1960s, voting was our organizing tool to demolish Jim Crow and achieve political impact. Since then, for me, it has been algebra. What’s math got to do with it? --- you ask. Everything, I say. 
\linebreak \linebreak 
Amidst the planetwide transformation we are undergoing, from industrial to information-age economies and culture, math performance has emerged as a critical measure of equal opportunity.}
\end{quotation}
\medskip
\medskip

Mathematical research also does not occur in a vacuum, especially when we apply our algorithms and other methods to social systems \cite{chiodo2023}. It is not a traditional part of curricular education or research education in the mathematical sciences to teach our students and the other junior community members about data ethics and other ethical considerations \cite{rycroft2022,muller2022}, but this needs to change.\footnote{One positive development is that some computer-science conferences now require authors to include ethics statements in their submitted papers.} What we choose to teach (and choose not to teach, or choose to mention in only a cursory way) impacts what our students and other mentees do with their education. Other disciplines (e.g., the social and medical sciences) have thought \emph{a lot} more about ethics than mathematics (and allied disciplines), and we should be guided by the best practices that they have developed. Many of these best practices arose in the aftermath of ethically problematic studies, and ethical guidelines have developed as people have tried to learn from past mistakes. As in other arenas, the mathematics community has had major ethical problems (e.g., through various forms of social toxicity), but these problems traditionally have arisen in issues other than the scientific use and abuse of data. This situation has changed; mathematicians are now also facing data-related issues directly in their research. What we learn, teach, and do needs to catch up to this reality.

It is important for mathematical scientists (and others) to think carefully about ethics and engage with it throughout their careers~\cite{chiodo2024}. We need to incorporate data ethics into the core education in the mathematical, statistical, and computational sciences. We are increasingly using social, animal, and human data (including potentially personal data) in our research. We need to think carefully about when it is appropriate to use such data and when it is not appropriate, and there needs to be systematic training to help mathematical scientists confront these issues. See \cite{loukides2018} for a short book on ethics and data science, \cite{ethics2016} for a theme journal issue on data ethics, \cite{veliz2021} for a wide-ranging handbook of digital ethics, \cite{hippo2021} for a discussion of the ethics landscape in mathematics and advocacy of a ``Hippocratic Oath'' for mathematics (also see \cite{hippo2022}), and \cite{skufca2021} for an article with a useful discussion and pointers to several relevant resources. See \cite{acm-ethics} for the Code of Ethics and Professional Conduct of the Association for Computing Machinery (ACM). Buell et al.~\cite{buell2022} leveraged the ethical-practice standards of the ACM and the American Statistical Association (ASA), which both represent disciplines that are relevant to mathematics, in a survey of mathematicians about ethical standards. See \cite{tractenberg2024a,tractenberg2024b} for further discussions of the guidelines that they developed for ethical mathematical practice, and see \cite{walk2024} for associated teaching materials. To help acknowledge issues like data ethics in the mathematical sciences, it is also relevant to include pertinent information on departmental websites. For an example website, which I helped produce at University of Oxford in the aftermath of a study \cite{Traud2012} by my collaborators and me that uses Facebook data, see \cite{MI-ox-data}.


\subsection{Course Material on Data Ethics and Society}
\label{courses}

Many existing courses discuss data ethics and the societal impact of data. Casey Fiesler (Department of Information Science, University of Colorado Boulder) has collected the syllabi of many such courses at \cite{fiesler2018}. A good book on ethics and data science to use as a starting point is \cite{loukides2018}, which the students in UCLA's course on societal impacts of data are asked to read.

Several courses have websites with much useful information about data ethics and related topics. I will highlight a few examples. Matt Salganik (Department of Sociology, Princeton University) has taught a graduate course (from which I drew some material for the present chapter) called \emph{Computational Social Science: Social Research in the Digital Age}. The material for the course's Fall 2016 edition is available at \cite{salganik2016}. It includes material that appeared later in book form \cite{salganik2017}. Johan Ugander (Department of Management Science and Engineering, Stanford University) has taught a graduate course called \emph{Data Ethics and Privacy}. See \cite{ugander2020} for the course website, from which I drew several articles that I mention in Section \ref{research}. Johan Ugander's graduate course \emph{Social Algorithms} \cite{Ugander2023} also has excellent resources. Chris Bail (Department of Sociology, Duke University) has taught a graduate course called \emph{Data Science \& Society}. See \cite{bail2022} for its course materials, schedule, and YouTube videos. A recent book by Chris Wiggins (Department of Applied Physics and Applied Mathematics, Columbia University) and Matthew L. Jones (Department of History, Columbia University) \cite{wiggins2023} is an adaptation of their course \emph{Data: Past, Present, and Future}~\cite{wiggins2024}.

Rachel Thomas (co-founder of {\sc fast.ai} and Professor of Practice in the Center for Data Science at Queensland University of Technology) has posted a wealth of resources at \cite{thomas-tweet2020}. This material includes a data-ethics course, a data-science blog, a diversity blog, and more. For example, take a look at her collection of short videos on ethics in machine learning \cite{thomas-videos}. The course \emph{Calling Bullshit} by Carl Bergstrom (Department of Biology) and Jevin West (Information School) at University of Washington includes a section on ethics \cite{bs2019}. This course also includes lots of other valuable material and reading suggestions.

The Cambridge University Ethics in Mathematics Project \cite{cambridge} also has relevant course material, including a description and recordings of an 8-lecture course 
called \emph{Ethics for the Working Mathematician}. Other useful material for courses on data ethics are research papers (see Section \ref{research}) and discussions of controversies (see Section \ref{controversy}).

\section{A Few Efforts Related to Data Ethics, Society, and Social Good}
\label{ucla}

In this section, I briefly discuss a few efforts that are related to data ethics and society. I highlight two recent efforts at University of California, Los Angeles (UCLA), and I briefly mention a few of the many efforts of others. I also highlight the importance of being careful and conscientious in these efforts.


\subsection{Two Recent Initiatives at UCLA}
\label{ucla-sub}
%


\subsubsection{An Undergraduate Course: Societal Impacts of Data}
\label{stats184}
I helped design a new undergraduate major in ``Data Theory''\footnote{See~\cite{datatheory} for a description of UCLA's undergraduate major in Data Theory.} at UCLA. All students in our Data Theory major must take a new upper-division course, which was designed by Mark Handcock of the Department of Statistics \& Data Science, called ``Societal Impacts of Data''. This course covers a variety of topics --- including privacy, algorithmic bias, and many others --- through an ethical lens. UCLA's course catalog has the following description of the course:
\begin{quotation}
\emph{Consideration of impacts that data collected today have upon individuals and society. Rapid increase in scale and types of data collected has impacted commerce and society in new ways. Consideration of economic, social and ethical, legal and political impacts of data, especially that collected on human behavior. Topics include privacy and data protection, intellectual property and confidentiality, sample selection and algorithms, equality and anti-discrimination.}
\end{quotation}
In a recent offering of the course, \cite{loukides2018} was used as reading material.

In my view, such courses should be mandatory not only for undergraduates who are majoring in data science (and similar topics), but also for all other students. Because data ethics is crucial for any member of human society, taking a course in it should be a core requirement for literally all students to obtain an undergraduate degree. However, it is especially important for the many students (in mathematics, statistics, computer science, and other subjects) who become data scientists or otherwise work with human data (and nonhuman animal data) in their careers.


\subsubsection{Social-Justice Data-Science Postdoctoral Scholars}
\label{sjds}

A new UCLA academic position, which I designed along with Deanna Needell (Department of Mathematics) and Mark Handcock, is a Social-Justice Data-Science (SJDS) postdoctoral scholar. This innovative postdoctoral position, which I hope to see in various forms at academic and other institutions worldwide, is a joint venture of UCLA's mathematics and statistics departments. An SJDS postdoc has two mentors: (1) a faculty member from the Department of Mathematics or the Department of Statistics \& Data Science and (2) a faculty member who is a social-justice scholar. We hope to continue to hire SJDS postdocs.

Mathematical, statistical, and computational scientists can accelerate the scientific study of social-justice issues by harnessing data. They can also leverage recent advances in data science into social justice and activism. Individuals who are trained in other fields (such as sociology) have a long tradition of such involvement in activism. Importantly, mathematicians, statisticians, and computer scientists also have a lot to contribute. Such contributions are not simply a matter of studying abstract problems in mathematical sociology or allied topics. Critically, it requires engagement with social-justice scholars and the communities and other stakeholders that we seek to help. This is what we expect SJDS postdocs to do. Moreover, by sponsoring SJDS scholars as postdocs in mathematics and statistics departments in positions of comparable prestige to their usual postdocs, the mathematical-science community (along with our colleagues in statistics, computer science, and other disciplines) can show students that these paths --- whether in academia, in industry, as a data scientist for a nonprofit organization that serves communities, or elsewhere --- are available and viable career pathways. It is crucial that we send this message.


\subsection{Some Other Noteworthy Efforts}
\label{other}

There are many other efforts (which take a variety of forms) on data and society. In Section \ref{courses}, I discussed several courses and resources on data ethics. In this subsection, I briefly highlight a few initiatives on data and society. 
 
 One of these efforts is Mechanism Design for Social Good (MD4SG) \cite{mech2018,md4sg}, which uses techniques from algorithms, optimization, mechanism design, and disciplinary insights to improve the equity, social welfare, and access to opportunities for historically underserved, disadvantaged, and marginalized communities. Another effort, which was launched by Timnit Gebru, is the Distributed AI Research Institute (DAIR) \cite{dair}. This institute works on community-based and interdisciplinary AI research. AI4All, which was founded by Fei-Fei Li and Olga Russakovsky in 2015, is a nonprofit organization that aims to increase diversity and inclusion in AI education, research, development, and policy \cite{AI4All}.
 
 
 \subsection{An Important Caveat}
\label{caveat}
 
There are a diverse variety of both research and practical efforts that fall under the auspices of ``data for social good'' \cite{abbasi2023}. Harnessing data for social good is both important and laudable, but there are many challenges \cite{abebe2020roles,buckee2022}. For example, it is vital to engage meaningfully and respectfully with communities, rather than attempting to help people by trying to be a ``new sheriff in town", which can be very harmful.


\section{Research Replication and Ethics, Open Data and Code, the Tension Between Privacy and Open Data and Code, and Some Controversy}
\label{open}


\subsection{Replication of Research}
\label{rep}

It is crucial to be able to replicate scientific research, and improvement in current practices is necessary to produce reproducible and reliable computational science \cite{coveney2021}. Obviously, it is necessary to be honest about data and other aspects of research, but professional obligations go far beyond mere honesty. It is important to precisely explain all details of analysis, implementation, and data cleaning in scholarly works; it is also important to openly provide source code and data. See \cite{abdill2024} for a set of recommendations for sharing code. The inclusion of precise explanations and source material is necessary for research dissemination. When a mathematician publishes a theorem and its proof, they also give other people a license to use it freely in their own work. The sharing and portability of knowledge lies at the core of both science and mathematics. Accordingly, it is a professional obligation to share research in a usable form, including by providing source code and data. Depriving readers of such material is analogous to publishing a theorem statement without a proof or publishing a theorem without permission to use it. 
 
To the extent possible, it is very important to publish the relevant data and code (including code to reproduce all figures) that accompanies a scholarly manuscript.\footnote{In the interest of admitting my own flaws, I note that I have been imperfect in my career about publishing source code with my papers. I am doing this increasingly often, and I seek to improve further.} In a manuscript, one also should explicitly and carefully discuss each step in the procedures for data anonymization, cleaning, sampling, and transformation. It is important to be explicit about anything that one does with data so that readers know precisely what choices have been made and can then evaluate whether or not they think that those choices are good ones for the analysis in a manuscript. For example, sampling biases can change the properties of data in fundamental ways \cite{stumpf2005}. Additionally, by providing the original data when possible, others can analyze that data in procedures that involve different choices. There are many choices that scientists make in data analysis --- it is impossible not to make such choices --- but these choices are a part of the scientific procedure in conducting research, so it is imperative to inform others of exactly what one has done (with particular highlighting of choices) in any scientific work. They may want to make different choices.

In making data publicly available, posting the output of synthetic models is safer than posting even the safest real-world data (see Section \ref{concern}). When using synthetic data, such as the output of numerical simulations of a differential equation or data that one generates from a random-network model, it is good to publish code to generate the output data (e.g., the examples) that is presented in a manuscript. Publishing user-friendly and open-source code is important for a paper's readers, and it also facilitates the fair evaluation of methods and results. One way to publish code is as supplementary material on a journal website; another way is through repositories such as Bitbucket, GitLab, and GitHub. Posting synthetic data is relevant not only for the output of numerical computations and any other data that one generates, but also even for examples (such as adjacency matrices in a paper about networks) that one constructs by hand in a definition--theorem--proof paper. How relevant it is to include such data depends on the sizes of the examples. Even posting the entries of a $10 \times 10$ matrix in a repository saves time for others and reduces transcription errors.

We are all human and it is easy to forget something or to inadvertently be insufficiently precise about a procedure, so gaps often occur. If somebody e-mails you to ask for a clarification, copy of code (even if poorly commented), or something else, it is important to respond and provide it to them (assuming that it is something that you have the legal and ethical right to provide). 


\subsection{Acknowledgements, Giving Credit, and Licenses to Use and Share Material}
\label{license}

Another part of doing scientific research and presenting it in scholarly works is acknowledging the contributions of others. Naturally, acknowledging contributions includes things like coauthorship and citations of prior work. It also includes things like acknowledging all sources of data, all sources of funding, and thanking people for their useful comments and ideas. In the acknowledgements section of a manuscript, one should include precise details of how one obtained data and how others can also obtain that data (especially if one cannot publish it oneself, perhaps because of privacy considerations or because it is not one's own data to share). It is important to be generous when acknowledging others in manuscripts. If somebody gives useful comments, one should thank them for it (assuming that they want to be thanked). 

One should be fair, appropriate, and precise when discussing prior work in a manuscript. It is crucial to give credit where it is due. The research in prior work has heterogeneous levels of mathematical rigor, scientific rigor, and even correctness. How one writes about such work is affected by such things, including some facets that are factual and others that may reflect a variety of opinions. For example, there is a difference in writing that something was ``shown'' versus ``reported'' in a prior work. The former wording has a built-in claim, by the writer of a manuscript, of the validity of that aspect of that prior work. By contrast, the latter is merely a historical fact (assuming that what one writes is itself accurate).

The increasing prominence of LLMs has brought new ethical concerns in acknowledgements and the creative process.\footnote{The use generative AI and LLMs also brings significant ethical concerns about data ownership and copyright infringement, such as through the training data that is used for AI-generated art and text \cite{francis2022,nyt2023}. It also can lead to comically disturbing situations, as illustrated by the recent discovery of child sexual-abuse material in a prominent data set (which has been downloaded by many researchers) in the AI community~\cite{csam2023} and by the AI-generated image of a rat with a gigantic penis in a (subsequently retracted) published scientific paper~\cite{rat2024}.} Scientific journals, funding agencies, professors, and universities are scrambling to develop policies that govern the use of LLMs and other text-generation tools in submissions \cite{brainard2023}. See \cite{MCM-AI} for the LLM and AI policy that the Consortium for Mathematics and its Applications (COMAP) has started using for the Mathematical Contest in Modeling (MCM) and their other contests. For UCLA's guidelines for the ethical and ``safe" use of AI, see UCLA's new AI website~\cite{ucla-AI}, which also discusses various AI resources.\footnote{Amidst the very serious concerns about AI, it is also important not to lose sight of its immense promise. For example, see the curated list \cite{AI-math} (and the announcement of it at \cite{tao-AI}) of resources for AI in mathematics (e.g., to use AI to assist in mathematical reasoning).}

Another aspect of open science is the different types of licensing that are available through Creative Commons. See \cite{creative} for a discussion of the different types of Creative Commons licenses. Some licenses allow work to be duplicated for any purpose, and other licenses enforce a variety of use restrictions. I advocate making one's work as open as possible, as others can then use it readily for purposes such as teaching and explaining ideas.


\subsection{Privacy Concerns and Practical Considerations for Open Data and Code}
\label{concern}

There are various tensions and practical considerations with the lofty ideals of openly publishing data and code. For example, one may not be allowed to publish data for privacy reasons or because of nondisclosure agreements. For empirical data, if you have permission to post something (e.g., does the data ``belong'' to somebody else?) and it does not pose privacy concerns, then it makes sense to post it because doing so promotes good science. Because of privacy issues, one may choose not to publish certain data that one is technically allowed to publish. It is crucial that researchers navigate these issues in a conscientious way. Another issue is that publicly posting usable code and data takes time and energy, and key participants (such as students and other junior researchers) in a project move on to other things, so the team members who are best equipped to do this effectively may no longer be available. In other words, there is sometimes a practical tension between publishing code and data and the well-being of one's junior collaborators. 

In Section \ref{rep}, I mentioned that it is important to indicate how one has anonymized data in scholarly works. When is data actually ``anonymous'', and is it ever possible to ``fully'' anonymize data? Consider the following scenario, which is discussed in \cite{salganik2016,salganik2017}. Suppose that we have a data set of medical records of individuals that includes their names, their home addresses, the zip codes of these addresses, their birth dates, their sexes, their ethnicities, the dates that each individual visited a doctor, the medical diagnoses, the medical procedures, and the prescribed medications. Now suppose that we ``anonymize'' this data set by removing the names and home addresses of all individuals. We now have a data set of ``anonymized'' medical records. Suppose that we obtain a data set of voting records and that this data set includes names, home addresses, political-party affiliations, voter-registration dates, zip codes, birth dates, and sexes. Both data sets include the zip codes, birth dates, and sexes of the individuals who are common to the two data sets. One can use such data --- namely, data that is common to these two data sets --- to de-anonymize people in the supposedly ``anonymized'' data set of medical records~\cite{arvind2008}. An infamous example of data de-anonymization by combining data sets led to the cancellation of the sequel to the Netflix Prize \cite{sequel}.

Given the simultaneous presence of privacy concerns and the desire to produce replicable scientific research, what should one do if the employed data, an algorithm (or part of an algorithm), or something else needs to remain private? This was one key topic of discussion following a publication by Bakshy et al.~\cite{bakshy2015}, who examined the exposure of different individuals to heterogeneous news and opinions in their Facebook feeds. The paper's authors, who were all Facebook employees at the time, could not reveal how Facebook determines the feeds that individuals see, so how can others replicate their work to try to evaluate and verify their observations and insights? Which of the insights in \cite{bakshy2015} apply exclusively to Facebook, and which of them also apply to other social-media platforms? In principle, it should be possible to attempt a weaker form of replication of the study's most interesting qualitative results, which are not merely a property of something that is specific to Facebook.


\subsection{Controversies in Studies and in Reactions to Studies}
\label{controversy}

Unsurprisingly, some studies that involve human data have been controversial. In other cases, there has been controversy in how authors of studies were treated in the aftermath of their work.

One controversial study, with much ensuing public discussion (see, e.g., \cite{luca2014} and many other sources), was an examination of emotional contagions using experiments with Facebook in which user feeds were altered~\cite{kramer2014}. There were angry accusations that the researchers manipulated people's emotions, with additional questioning of subsequent actions by the journal that published the paper. There were also discussions of the procedure to obtain permission to undertake the study in the first place. One key consideration is that there are crucial differences between the ethical procedures for academic and commercial researchers \cite{grindrod2016}. Academic researchers need to obtain approval from an Institutional Review Board (see Section \ref{irb}) before undertaking a study like this that involves humans, whereas companies like Facebook have publication review boards to approve the publication of a study after it has already been done. Therefore, we know that this study occurred because Facebook concluded that it could be published. By contrast, we do not know about what research is done with our data by Facebook and other entities when an associated document is not placed in the public domain. This leads to an important question: Should academic researchers and companies follow the same rules?

Many technology companies have units that do research on data ethics and related subjects, although that too can lead to controversy. One example is the departure from Google of researcher Timnit Gebru and others who study data ethics \cite{gebru-google} in the aftermath of a paper that Gebru and her collaborators wrote about the significant risks (including environmental costs, unknown and dangerous biases, and potential uses to deceive people) of LLMs~\cite{hao2020}.

Other research about online social networks has also been controversial. I will mention two examples that have been influential scientifically because of their research findings. One of these examples involved experimental manipulation of feeds by seeding posts on social media with a small number of fake initial upvotes or fake initial downvotes \cite{aral2012}. In this study, the researchers found that initial upvotes had a persistent effect on the overall positivity of the votes of posts, whereas the initial downvotes were overturned. The other example is the ``Tastes, Ties, and Times'' study of several waves of Facebook data from students at an Ivy League university in the United States \cite{ttt2008}. See \cite{ttt-news} for one discussion of the data-privacy controversy of this study and its associated data set.


\subsection{Summary}
\label{open-summary}

To summarize some key ideas in Section \ref{open}, here are a few things to think about:
\begin{itemize}
\item{There is tension between open data and personal privacy.}
\item{The use and publication of data, code, and anything else that one reports in a manuscript or discusses with others can be subject to terms-of-service agreements and nondisclosure agreements.}
\item{In what sense can you make your research replicable if you cannot make all of your data (or algorithms or something else) publicly available? There are weaker notions of replication, such as whether or not others can observe similar phenomena in circumstances that are similar but not precisely the same. For instance, in studying human behavior on social media, perhaps certain phenomena are very similar on Facebook and 
$\mathbb{X}$ (which formerly was called Twitter), but other phenomena are specific to only one of these two social-media platforms.}
\item{Are you comfortable doing research in collaboration with private companies or government entities? Maybe there are some entities with which you are willing to collaborate (perhaps depending on their purposes, goals, and history), but there are others with which you are not willing to collaborate or use data from? If you work with or for such an entity, what is permissible to include in a publication or post online?}
\end{itemize}


\section{Ethical Principles and Other Considerations in the Scientific Use of Human Data}
\label{principle}

In this section, I discuss some ethical principles and other considerations in the scientific use of human data.\footnote{Parts of my discussion draw heavily from material in Matt Salganik's course on computational social science \cite{salganik2016}. See his associated book \cite{salganik2017} for further discussion.} I also discuss Institutional Review Boards (IRBs).\footnote{At some institutions, the backronym for IRB is Internal Review Board.}
Much of my discussion also applies to data from nonhuman animals, but I am focusing on human data in this chapter, so I typically phrase my exposition in human terms. 


\subsection{Online Courses for Ethics Training}\label{training}

In studies that use human data, it is important to think carefully about ethics and to have formal training in it. A popular choice is the Collaborative Institutional Training Initiative (CITI) program \cite{citi}, which offers a variety of courses. For more information, see the website of the UCLA Office of the Human Research Protection Program (OHRPP) \cite{ohrpp}. See Section \ref{irb} for further discussion of these courses and of IRBs.


\subsection{Four Key Principles}\label{key}

I now enumerate a few critical ethical principles, which are discussed in detail in chapter 6 of \cite{salganik2017}.

In scientific pursuits that involve human data, it is important to do the following:
\begin{itemize}
\item{be honest and fair (obviously);}
\item{design ethically thoughtful research;}
\item{explain your decisions to others.}
\end{itemize}
Four key principles in research that involves humans are
\begin{itemize}
\item{respect for persons;}
\item{beneficence;}
\item{justice;}
\item{respect for law and public interest.}
\end{itemize}
These four principles can come into tension with each other, so how do we balance them?

In conducting research with human data, there is a sliding scale: the more your research has the potential to violate personal privacy, the more helpful for humanity its outcome needs to have the potential to be. Four things to ponder with research that involves personal data are the following:
\begin{itemize}
\item{informed consent;}
\item{understanding and managing informational risk;}
\item{privacy;}
\item{making decisions in the face of uncertainty.}
\end{itemize}
As you design and conduct research, put yourself in the shoes of other people. Think of research ethics as continuous (i.e., there is a sliding scale), rather than as discrete.


\subsection{Institutional Review Boards (IRBs)}\label{irb}

\begin{quotation}
``Yeah, yeah, but your scientists were so preoccupied over whether or not they could, they didn't stop to think if they \emph{should}."  \\ \\ 
(stated by the character Ian Malcolm in \emph{Jurassic Park})
\end{quotation}

\medskip

In many situations, it is a professional requirement to obtain permission to undertake a study in the first place. Such permission gives an ethical floor to satisfy; it is not a ceiling. It may be legally and professionally permissible to do something, but it is important to hold oneself to higher standards when it comes to whether or not it is actually the right decision to do it. For example, the limits to informed consent with human data~\cite{lovato2022} may influence such a decision. Whatever you decide for your own work, make sure that you think carefully about it.

In universities in the United States, a researcher who is working with personal data needs to check with their university's IRB to ensure that they are conducting research in an ethical way. Universities in other countries, private companies, government laboratories, and other organizations often have bodies that are similar to IRBs, but the procedures and especially the specific details are very different \cite{grindrod2016}. A university IRB may inform you that you do not need to submit a formal application for a research project to be approved, or they may inform you that a formal application is necessary. Let your IRB know what data you have (or what data you plan to acquire and how you plan to acquire it) and what you plan to do with it. Different IRBs can rule differently. When an IRB grants permission to undertake a study, they have decided that a proposed project is above the ethical floor and hence that it is permissible to do that project. One's own standards should be higher. (Again, ethical approval is a floor, rather than a ceiling.) In this light, it is worth examining the discussion following a controversial IRB-approved study of ``emotional manipulation'' in which researchers adjusted user feeds on Facebook \cite{luca2014,kramer2014}.

For a more concrete idea about IRBs and conducting ethical research in a university setting, see the online materials at UCLA's Office of the Human Research Protection Program (OHRPP) \cite{ohrpp}. The requirements to conduct research with human data (and nonhuman animal data) include taking various online training courses, such as those in the CITI program. These courses, which are available online at \cite{citi}, are in common use in the United States. The training that is required, expected, and available for research projects that involve human data and other sensitive data differs substantially in different countries. Some of the topics that are covered in courses on ethical research are animal care and use, biosafety and biosecurity, human-subject research, information privacy and security, and responsible conduct of research. It is useful to take a variety of these courses even when they are not required.


\subsection{Another Salient Warning}
\label{salient}

It is important to be cognizant that your research can potentially be ``weaponized'' by other people --- who may not care about nuances in research findings and who may interpret unfortunate wording or insufficiently careful exposition in nefarious ways --- so you need to be conscientious about the precise wording in your publications and other media. This is particularly relevant in research on social systems and on tools that are applicable to social systems~\cite{stemwedel2022}. This possibility also makes it particularly crucial to be open-minded about potential sources of bias in your research. As always, it is important to be ethically thoughtful.


\section{Some Research and Lay Articles that are Relevant to Data Ethics and Data Privacy}
\label{research}

There is a wealth of research about data ethics and related topics in computer science, sociology, and other fields \cite{jobin2019,lazer2021,wagner2021,abebe2020roles}. There is no way that I can possibly be exhaustive,\footnote{I selected some of the highlighted articles, including both research studies and lay articles, from the website for Johan Ugander's course on data ethics and data privacy \cite{ugander2020}.} so I will highlight a selection of studies to give a taste of existing research. I will also briefly discuss a few articles in blogs, news websites, and other nontechnical venues.
Unsurprisingly, much of this research is simultaneously fascinating and concerning \cite{xkcd2072}.


\subsection{Algorithmic Biases and Other Biases}
\label{biases}

A key research area is the analysis and mitigation of biases in computational techniques and data analysis. 

One central topic is algorithmic bias, in which systematic and repeatable errors result in unfair outcomes, such as privileging one group of people over others \cite{birhane2021}. For example, such biases have systematically hurt certain racial and demographic groups in predictive policing \cite{proactive2018}. Algorithmic biases also reinforce negative racial and gender biases in online search-engine results~\cite{noble2018}. Another example is an algorithm for prioritizing COVID-19 vaccinations that severely disadvantaged human-facing medical workers \cite{covid2020}. 

There is much research on the mitigation algorithmic bias, such as in the algorithmic hiring of people for jobs \cite{raghavan2020} and in classifying individuals (e.g., for admission to a university) \cite{dwork2012}. To encourage transparent reporting, clarify intended use cases, and minimize use in inappropriate contexts, some researchers have proposed the inclusion of ``model cards'' to accompany trained machine-learning models~\cite{mitchell2019-model-cards}. Such model cards are short documents that give intended use cases; benchmarked evaluations across cultural, demographic, phenotypic, and intersectional groups; and other salient information. 

Biases in computer systems, machine learning, and AI go far beyond only algorithmic bias \cite{friedman1996,hogg2024}. See \cite{mehrabi2021,barocas2023,corbett2023} for reviews of bias and ``fairness" in machine learning, and see \cite{ntoutsi2020} for an introductory survey of bias in data-driven AI systems. As advocated in \cite{kasy2021}, it is important to go beyond notions of mere algorithmic ``fairness'' (which focuses on intra-group versus inter-group differences). One must also analyze (1) inequality and the causal impact of algorithms and (2) the distribution of power. It is also important to be cognizant of traps that can beset work on ``fair'' machine learning in sociotechnical systems \cite{selbst2019}. In ``fair" machine learning, it is necessary to carefully examine context-specific consequences \cite{corbett2023}.

To mitigate biases, it is imperative for researchers to carefully justify their choices of data sets. It is crucial to consider the social context (country, gender, race, and so on) of data. Why is one is using a particular data set, and is it appropriately representative for the problem under study \cite{koch2021}? It is common and often convenient to import data sets that were used originally to study one problem for investigations of many other problems, and that can lead to several issues. This is particularly relevant for human data and social data, and it is important to think carefully about whether a data set is appropriate for a given study. Particular facets of a setting can interfere with attempts to generalize a study's insights from its specific use case to other potentially similar situations.


\subsection{Human Privacy, Personal Characteristics, and Personalization}
\label{privacy-char}

Other salient research focuses on human privacy. As has been studied thoroughly \cite{engle2016}, people are tracked extensively when they visit websites. User profiles, which encode both characteristics and behavior, are generated through interactions with websites and other digital systems. User profiling is ubiquitous in daily life, and it is also a vast area of research~\cite{purificato2024}. It is possible to infer private traits and attributes from digital records of human behavior \cite{kosinski2013}. One can also steal the identities of visitors to websites \cite{arvind2010}, and trackers can use social-network structure to de-anonymize data from browsing the World Wide Web \cite{su2017}. Companies can exploit the information that they observe (or infer) from website visits. For example, about a decade ago, the online travel agent Orbitz showed higher prices for flights to users of Macintosh computers than to users of other types of computers \cite{white2012}. It is imperative that individuals and other stakeholders have agency in how their personal data are used \cite{sadowski2021}, although this can come into tension with the scholarly desire to promote open data (see Section \ref{open}). For example, the sharing of data from the continent of Africa has often been driven by non-African stakeholders \cite{abebe2021}. Additionally, conventional studies of algorithmic fairness are centered on Western culture; data proxies are different in different cultures, and it is important to consider local context when building data models \cite{samba2021}. See \cite{lazer2021,wagner2021} for discussions of access, ethics, and best practices in the algorithmic measurement of human data.

There are many ways to infer the characteristics of people and communities, as well as social and other ties between people, using data analysis. For example, it is possible to infer social ties by examining the geographical proximity of individuals in time and space using offline or online data \cite{crandall2010}. Additionally, researchers have used machine learning and Google Street View to estimate the demographic composition of neighborhoods across the United States \cite{gebru2017}. Unsurprisingly, there are significant gaps and biases (and associated exclusion issues) in such geographic data \cite{graham2022}. Moreover, as I discussed in Section \ref{open}, the fact that individuals appear in multiple social networks (e.g., multiple social-media platforms and multiple databases) gives considerable ability to identify them or information about them even when they appear to be ``anonymous'' in those individual networks \cite{arvind2008}. For example, the inference of the sexual orientation of some individuals by combining Netflix data --- which was released as part of a public competition to improve Netflix's algorithmic ability to infer user ratings of movies --- of movie rentals with movie rating data from the Internet Movie Database (IMDB) led to legal troubles, and Netflix ultimately cancelled a planned sequel competition \cite{sequel}. More recently, researchers have successfully inferred gender using mobile-payment data \cite{stobaugh2023}. Moreover, one can use interactions to identify people even across long time periods \cite{cretu2022}. Tools like data analysis and machine learning can yield crucial and actionable insights (e.g., of racial disparities in police stops \cite{pierson2020}), but they need to be applied in a careful and respectful way.

A related issue is personalization \cite{jugander2017medium}, such as in the targeted advertisements that inundate people on social-media platforms. These advertisements can inform us of desirable products (such as t-shirts, dice, or plushies in my case), but the online personalization of commercial and political advertising can cause serious problems. Infamously, Donald Trump's 2016 presidential campaign employed the company Cambridge Analytica for targeted online political advertising. Additionally, the movie \emph{Straight Outta Compton} was advertised differently to people based on their inferred demographic characteristics \cite{compton}. There is much pertinent research. For example, there are studies of the effectiveness of spam-based marketing \cite{spam2009}, the detection of spam in social-bookmarking websites (such as Pinterest and Digg)~\cite{markines2009}, and other aspects of spam and related online marketing. Users on social-media platforms like Facebook have personalized feeds, and a study that involved the adjustment of such feeds \cite{kramer2014} led to a major controversy, including with the issue of manipulating people's emotions \cite{luca2014}.


\subsection{AI Ethics Guidelines, Tradeoffs in Different Ethical Values, and Theoretical Barriers}
\label{guidelines}

Guidelines for AI ethics appear to have converged on five key values:\footnote{These values overlap with but are organized differently from the five principles to guide the design, use, and deployment of automated systems in the recent blueprint for an AI Bill of Rights from the United States government's Office of Science and Technology Policy (OSTP)~\cite{OSTP2022}. The OSTP's five principles are (1) safe and effective systems, (2) algorithmic-discrimination protections, (3) data privacy, (4) notice and explanation, and (5) human alternatives, consideration, and fallback.} transparency, fairness, safety, accountability, and privacy~\cite{jobin2019}. There are inevitable tensions between these ethical values. Academics, public-sector organizations, and private companies emphasize different ethical values, with particularly systematic differences between practitioners and the general public~\cite{jakesch2022}. 
 
 It is mathematically impossible for an AI system to simultaneously include many parameters, be robust to poisoning (e.g., with fake data) by an adversarial actor, and preserve privacy~\cite{mahdi2023}, so the use of AI will always include tradeoffs between different ethical values even when everything works perfectly. Of course, AI systems don't work perfectly. For example, it was proven recently that one can plant undetectable backdoors into machine-learning models and thereby confuse ``adversarially robust" classifiers, demonstrating a major theoretical barrier to the certification of adversarial robustness~\cite{goldwasser2022}.


\section{Conclusion}
\label{summary}

The use of digital data to study and analyze humans, create technologies, and enact policies that involve humans is both powerful and dangerous \cite{editorial2021,OSTP2022}. Mathematical scientists have a younger tradition of studying human data than researchers in many other disciplines, such as the social and medical sciences. Most researchers in the mathematical, statistical, and computational sciences have not had research-ethics and data-ethics training. This needs to change. One recent paper has even proposed a ``Hippocratic Oath'' for the mathematical sciences \cite{hippo2021}. The mathematical-science community needs to learn from the best practices of other disciplines to ensure that our research is ethical. Other scholarly communities have been considering research ethics a lot longer than we have, and it is important that we learn from them. As in the social and medical sciences, the mathematical and computational sciences need a robust program of ethics training.

In this chapter, I have discussed many aspects of data ethics in education, publishing, and research. In some of my discussions, I have also touched on issues of data and justice. It is important to distinguish between ethics and justice, including when it comes to data~\cite{kitchin2014}. Ethics presupposes laws and social norms, and then one considers what is morally appropriate within those frameworks. Justice is concerned with how to change laws and modify social norms to attain broader equity. Both ethics and justice are core aspects of our interactions with data.

I encourage you to read widely, think about, and discuss how to do research ethically, especially for studies of --- or with consequences for --- social systems and human data. It is also valuable to read about past research and societal controversies. There have been mistakes in the past (and there continue to be mistakes), and we need to learn from them. We may all set our ethical bars in different places and have different views on different issues, but our scholarship needs to be conscientious and ethically thoughtful. As a reminder, official approval (e.g., from an IRB) to undertake a study is only a lower bound. The ethical bar to surpass in research design and performance is a sliding bar: the more potential for invasion of human privacy (or other potential harm), the more valuable to humanity the potential outcome of a research project has to be.

Be ethically thoughtful.



\begin{thebibliography}{WSAO{\etalchar{+}}21}

\bibitem[AAB{\etalchar{+}}21]{abebe2021}
Rediet Abebe, Kehinde Aruleba, Abeba Birhane, Sara Kingsley, George Obaido,
  Sekou~L Remy, and Swathi Sadagopan, \emph{Narratives and counternarratives on
  data sharing in {A}frica}, Proceedings of the 2021 ACM Conference on
  Fairness, Accountability, and Transparency (FAccT '21) (New York, NY, USA),
  Association for Computing Machinery, 2021, pp.~329--341.

\bibitem[ABK{\etalchar{+}}20]{abebe2020roles}
Rediet Abebe, Solon Barocas, Jon Kleinberg, Karen Levy, Manish Raghavan, and
  David~G. Robinson, \emph{Roles for computing in social change}, Proceedings
  of the 2020 Conference on Fairness, Accountability, and Transparency (New
  York City, NY, USA), FAT* '20, Association for Computing Machinery, 2020,
  pp.~252--260.

\bibitem[ACX23]{abbasi2023}
Ahmed Abbasi, Roger H.~L. Chiang, and Jennifer Xu, \emph{Data science for
  social good}, Journal of the Association for Information Systems \textbf{24}
  (2023), no.~6, 1439--1458.

\bibitem[AG18]{mech2018}
Rediet Abebe and Kira Goldner, \emph{Mechanism design for social good}, AI
  Matters \textbf{4} (2018), no.~3, 27--34.

\bibitem[AI4]{AI4All}
\emph{{AI4All}}, available at \url{https://ai-4-all.org} (accessed 26 December
  2023).

\bibitem[{Ass}18]{acm-ethics}
{Association for Computing Machinery}, \emph{{ACM Code of Ethics and
  Professional Conduct}}, 2018, available at \url{https://ethics.acm.org}
  (accessed 26 December 2023).

\bibitem[ATG24]{abdill2024}
Richard~J. Abdill, Emma Talarico, and Laura Grieneisen, \emph{A how-to guide
  for code-sharing in biology}, arXiv:2401.03068 (2024).

\bibitem[AW12]{aral2012}
Sinan Aral and Dylan Walker, \emph{Identifying influential and susceptible
  members of social networks}, Science \textbf{337} (2012), no.~6092, 337--341.

\bibitem[Bai22]{bail2022}
Christopher~A. Bail, \emph{{Data Science \& Society}}, 2022, Sociology 367S,
  Spring 2022, Duke University. Available at \url{https://dssoc.github.io}.

\bibitem[BBS22]{buckee2022}
Caroline Buckee, Satchit Balsari, and Andrew Schroeder, \emph{Making data for
  good better}, PLoS Digital Health \textbf{1} (2022), no.~1, e0000010.

\bibitem[BFPV21]{short2021}
Heather~Z. Brooks, Michelle Feng, Mason~A. Porter, and Alexandria Volkening,
  \emph{{Short Course: Mathematical and Computational Methods for Complex
  Social Systems}}, 2021, available at
  \url{https://zerodivzero.com/short_course/aaac8c66007a4d23a7aa14857a3b778c/titles}.

\bibitem[BGMMS21]{parrots2021}
Emily~M. Bender, Timnit Gebru, Angelina McMillan-Major, and Shmargaret
  Shmitchell, \emph{On the dangers of stochastic parrots: {C}an language models
  be too big?~\parrot}, Proceedings of the 2021 ACM Conference on Fairness,
  Accountability, and Transparency (New York City, NY, USA), FAccT '21,
  Association for Computing Machinery, 2021, pp.~610--623.

\bibitem[BHN23]{barocas2023}
Solon Barocas, Moritz Hardt, and Arvind Narayanan, \emph{{Fairness and Machine
  Learning: Limitations and Opportunities}}, MIT Press, Cambridge, MA, USA,
  2023, available at \url{http://www.fairmlbook.org}.

\bibitem[Bir21]{birhane2021}
Abeba Birhane, \emph{Algorithmic injustice: {A} relational ethics approach},
  Patterns \textbf{2} (2021), no.~2, 100205.

\bibitem[BMA15]{bakshy2015}
Eytan Bakshy, Solomon Messing, and Lada~A Adamic, \emph{Exposure to
  ideologically diverse news and opinion on {Facebook}}, Science \textbf{348}
  (2015), no.~6239, 1130--1132.

\bibitem[BPT22]{buell2022}
Catherine~A. Buell, Victor~I. Piercey, and Rochelle~E. Tractenberg,
  \emph{Leveraging guidelines for ethical practice of statistics and computing
  to develop standards for ethical mathematical practice: {A} white paper},
  arXiv:2209.09311 (2022).

\bibitem[Bra23]{brainard2023}
Jeffrey Brainard, \emph{As scientists explore {AI}-written text, journals
  hammer out policies}, Science (2023), 22 February, available at
  \url{https://www.science.org/content/article/scientists-explore-ai-written-text-journals-hammer-policies}.

\bibitem[BW19]{bs2019}
Carl~T. Bergstrom and Jevin West, \emph{{Calling Bullshit: Data Reasoning in a
  Digital World. Week 10. The Ethics of Calling Bullshit.}}, 2019, Informatics
  270 / Biology 270, University of Washington. Available at
  \url{https://www.callingbullshit.org/syllabus.html#Ethics}.

\bibitem[cam]{cambridge}
\emph{{The Cambridge University Ethics in Mathematics Project}}, available at
  \url{https://www.ethics.maths.cam.ac.uk} (accessed 26 December 2023).

\bibitem[CBC{\etalchar{+}}10]{crandall2010}
David~J. Crandall, Lars Backstrom, Dan Cosley, Siddharth Suri, Daniel
  Huttenlocher, and Jon Kleinberg, \emph{Inferring social ties from geographic
  coincidences}, Proceedings of the National Academy of Sciences of the United
  States of America \textbf{107} (2010), no.~52, 22436--22441.

\bibitem[CDNSG23]{corbett2023}
Sam Corbett-Davies, Hamed Nilforoshan, Ravi Shroff, and Sharad Goel, \emph{The
  measure and mismeasure of fairness}, Journal of Machine Learning Research
  \textbf{24} (2023), 312.

\bibitem[CGH21]{coveney2021}
Peter~V. Coveney, Derek Groen, and Afons~G. Hoekstra, \emph{Reliability and
  reproducibility in computational science: {I}mplementing validation,
  verification and uncertainty quantification \emph{in silico}}, Philosophical
  Transactions of the Royal Society A \textbf{379} (2021), 20200409.

\bibitem[Che20]{covid2020}
Caroline Chen, \emph{Only seven of {S}tanford's first 5,000 vaccines were
  designated for medical residents}, ProPublica (2020), available at
  \url{https://www.propublica.org/article/only-seven-of-stanfords-first-5-000-vaccines-were-designated-for-medical-residents}.

\bibitem[{CIT}]{citi}
{CITI Program}, \emph{{Explore Our Courses}}, available at
  \url{https://about.citiprogram.org} (accessed 26 December 2023).

\bibitem[CM23]{chiodo2023}
Maurice Chiodo and Dennis M\"{u}ller, \emph{Manifesto for the responsible
  development of mathematical works --- {A} tool for practitioners and for
  management}, arXiv:2306.09131 (2023).

\bibitem[CM24]{chiodo2024}
\bysame, \emph{A field guide to ethics in mathematics}, Notices of the American
  Mathematical Society \textbf{71} (2024), no.~7, 939--947.

\bibitem[CMM{\etalchar{+}}22]{cretu2022}
Ana-Maria Crețu, Federico Monti, Stefano Marrone, Xiaowen Dong, Michael
  Bronstein, and Yves-Alexandre {de Montjoye}, \emph{Interaction data are
  identifiable even across long periods of time}, Nature Communications
  \textbf{3} (2022), 313.

\bibitem[Col23]{csam2023}
Samantha Cole, \emph{Largest dataset powering {AI} images removed after
  discovery of child sexual abuse material}, 404 Media (2023), 20 December,
  available at
  \url{https://www.404media.co/laion-datasets-removed-stanford-csam-child-abuse/}.

\bibitem[{Con}]{MCM-AI}
{Consortium for Mathematics and its Applications}, \emph{Use of large language
  models and generative {AI} tools in {COMAP} contests}, available at
  \url{https://www.contest.comap.com/undergraduate/contests/mcm/flyer/Contest_AI_Policy.pdf}
  (accessed 21 April 2024).

\bibitem[D'A22]{hippo2022}
Susan D'Agostino, \emph{Do scientists need an {AI} {H}ippocratic oath? {M}aybe.
  {M}aybe not.}, Bulletin of the Atomic Scientists (2022), available at
  \url{https://thebulletin.org/2022/06/do-scientists-need-an-ai-hippocratic-oath-maybe-maybe-not/}.

\bibitem[dai]{dair}
\emph{{Distributed AI Research Institute (DAIR)}}, available at
  \url{https://www.dair-institute.org} (accessed 26 December 2023).

\bibitem[DHP{\etalchar{+}}12]{dwork2012}
Cynthia Dwork, Moritz Hardt, Toniann Pitassi, Omer Reingold, and Richard Zemel,
  \emph{Fairness through awareness}, Proceedings of the 3rd Innovations in
  Theoretical Computer Science Conference (New York City, NY, USA), ITCS '12,
  Association for Computing Machinery, 2012, pp.~214--226.

\bibitem[Dor20]{precious}
Lori Dorn, \emph{A deepfaked {Gollum} performs a `precious' cover of {Sinead
  O'Connor's} version of {`Nothing Compares 2 U'}}, Laughing Squid (2020),
  available at
  \url{https://laughingsquid.com/gollum-sinead-oconnor-nothing-compares-2-u/}
  (accessed 26 December 2023). The music video itself is available at
  \url{https://www.youtube.com/watch?v=9d9pZi7uZQs} (accessed 26 December
  2023).

\bibitem[Edi21]{editorial2021}
Editors, \emph{The powers and perils of using digital data to understand human
  behaviour}, Nature \textbf{595} (2021), 149--150.

\bibitem[EMFG{\etalchar{+}}23]{mahdi2023}
El-Mahdi El-Mhamdi, Sadegh Farhadkhani, Rachid Guerraoui, Nirupam Gupta,
  L\^{e}-Nguy\^{e}n Hoang, Rafael Pinot, S\'{e}bastien Rouault, and John
  Stephan, \emph{On the impossible safety of large {AI} models},
  arXiv:2209.15259 (2023).

\bibitem[EN16]{engle2016}
Steven Englehardt and Arvind Narayanan, \emph{Online tracking: {A}
  1-million-site measurement and analysis}, Proceedings of the 2016 ACM SIGSAC
  Conference on Computer and Communications Security (New York, NY, USA), CCS
  '16, Association for Computing Machinery, 2016, pp.~1388--1401.

\bibitem[ET19]{engin2019}
Zeynep Engin and Philip Treleaven, \emph{Algorithmic government: {A}utomating
  public services and supporting civil servants in using data science
  technologies}, The Computer Journal \textbf{62} (2019), no.~3, 448--460.

\bibitem[eth]{ethicalmath}
\emph{{Ethical Mathematics}}, available at \url{https://ethicalmath.com}
  (accessed 26 December 2023).

\bibitem[{Exp}18]{xkcd2072}
{Explain xkcd}, \emph{{2072: Evaluating Tech Things}}, 2018, available at
  \url{https://www.explainxkcd.com/wiki/index.php/2072:_Evaluating_Tech_Things}.
  The original comic (by Randall Munroe) is available at
  \url{https://xkcd.com/2072/}.

\bibitem[Fer24]{ferrara2023}
Emilio Ferrara, \emph{{GenAI} against humanity: {N}efarious applications of
  generative artificial intelligence and large language models}, Journal of
  Computational Social Science (2024), available at
  \url{https://doi.org/10.1007/s42001-024-00250-1}.

\bibitem[Fie18]{fiesler2018}
Casey Fiesler, \emph{Tech ethics curricula: {A} collection of syllabi}, Medium
  (2018), available at
  \url{https://cfiesler.medium.com/tech-ethics-curricula-a-collection-of-syllabi-3eedfb76be18}.

\bibitem[FN96]{friedman1996}
Batya Friedman and Helen Nissenbaum, \emph{Bias in computer systems}, ACM
  Transactions on Information Systems \textbf{14} (1996), no.~3, 330--347.

\bibitem[fPM20]{ipam-fakery}
Institute for Pure and Applied Mathematics, \emph{{White Paper: Deep Fakery}},
  2020, available at
  \url{http://www.ipam.ucla.edu/news/white-paper-deep-fakery/}.

\bibitem[Fra22]{francis2022}
Matthew~R. Francis, \emph{The ethics of artificial intelligence-generated art},
  SIAM News (2022), November, available at
  \url{https://sinews.siam.org/Details-Page/the-ethics-of-artificial-intelligence-generated-art}.

\bibitem[FT16]{ethics2016}
L.~Floridi and M.~Taddeo, \emph{What is data ethics?}, Philosophical
  Transactions of the Royal Society A \textbf{374} (2016), no.~2083, 20160360.

\bibitem[GD22]{graham2022}
Mark Graham and Martin Dittus, \emph{{Geographies of Digital Exclusion: Data
  and Inequality}}, Pluto Press, London, UK, 2022.

\bibitem[GKVZ22]{goldwasser2022}
Shafi Goldwasser, Michael~P. Kim, Vinod Vaikuntanathan, and Or~Zamir,
  \emph{Planting undetectable backdoors in machine learning models}, 2022 IEEE
  63rd Annual Symposium on Foundations of Computer Science (FOCS) (Los
  Alamitos, CA, USA), IEEE Computer Society, 2022, pp.~931--942.

\bibitem[GKW{\etalchar{+}}17]{gebru2017}
Timnit Gebru, Jonathan Krause, Yilun Wang, Duyun Chen, Jia Deng, Erez~Lieberman
  Aiden, and Li~Fei-Fei, \emph{{Using deep learning and {Google Street View} to
  estimate the demographic makeup of neighborhoods across the {United
  States}}}, Proceedings of the National Academy of Sciences of the United
  States of America \textbf{114} (2017), no.~50, 13108--13113.

\bibitem[GM23]{nyt2023}
Michael~M. Grynbaum and Ryan Mac, \emph{{\emph{The Times}} sues {OpenAI} and
  {Microsoft} over {A.I.} use of copyrighted work}, The New York Times (2023),
  27 December, available at
  \url{https://www.nytimes.com/2023/12/27/business/media/new-york-times-open-ai-microsoft-lawsuit.html}.

\bibitem[Gri16]{grindrod2016}
Peter Grindrod, \emph{Beyond privacy and exposure: {E}thical issues within
  citizen-facing analytics}, Philosophical Transactions of the Royal Society of
  London A: Mathematical, Physical and Engineering Sciences \textbf{374}
  (2016), no.~2083, 20160132.

\bibitem[Hao20]{hao2020}
Karen Hao, \emph{We read the paper that forced {Timnit Gebru out of Google}.
  {H}ere's what it says.}, MIT Technology Review (2020), available at
  \url{https://www.technologyreview.com/2020/12/04/1013294/google-ai-ethics-research-paper-forced-out-timnit-gebru/}.

\bibitem[Has21]{hassel2021}
Gry Hasselbalch, \emph{{Data Ethics of Power: A Human Approach in the Big Data
  and AI Era}}, Edward Elgar Publishing, Northampton, MA, USA, 2021.

\bibitem[HV24]{hogg2024}
David Hogg and Soledad Villar, \emph{{\it Position:} {I}s machine learning good
  or bad for the natural sciences?}, arXiv:2405.18095 (2024).

\bibitem[JBAO22]{jakesch2022}
Maurice Jakesch, Zana Bu\c{c}inca, Saleema Amershi, and Alexandra Olteanu,
  \emph{How different groups prioritize ethical values for responsible {AI}},
  Proceedings of the 2022 ACM Conference on Fairness, Accountability, and
  Transparency (New York City, NY, USA), FAccT '22, Association for Computing
  Machinery, 2022, pp.~310--323.

\bibitem[JIV19]{jobin2019}
Anna Jobin, Marcello Ienca, and Effy Vayena, \emph{The global landscape of {AI}
  ethics guidelines}, Nature Machine Intelligence \textbf{1} (2019), 389--399.

\bibitem[KA21]{kasy2021}
Maximilan Kasy and Rediet Abebe, \emph{Fairness, equality, and power in
  algorithmic decision-making}, Proceedings of the 2021 ACM Conference on
  Fairness, Accountability, and Transparency (FAccT '21) (New York, NY, USA),
  Association for Computing Machinery, 2021, pp.~576--586.

\bibitem[KDHF21]{koch2021}
Bernard Koch, Emily Denton, Alex Hanna, and Jacob~G. Foster, \emph{Reduced,
  reused and recycled: {T}he life of a dataset in machine learning research},
  NeurIPS 2021: Proceedings of the Neural Information Processing Systems Track
  on Datasets and Benchmarks (J.~Vanschoren and S.~Yeung, eds.), 2021.

\bibitem[KGH14]{kramer2014}
Adam D.~I. Kramer, Jamie~E. Guillory, and Jeffrey~T. Hancock,
  \emph{Experimental evidence of massive-scale emotional contagion through
  social networks}, Proceedings of the National Academy of Sciences of the
  United States of America \textbf{111} (2014), no.~24, 8788--8790.

\bibitem[Kit14]{kitchin2014}
Rob Kitchin, \emph{{The Data Revolution: Big Data, Open Data, Data
  Infrastructures \& Their Consequences}}, Sage Publishing, Thousand Oaks, CA,
  USA, 2014.

\bibitem[KKL{\etalchar{+}}09]{spam2009}
Chris Kanich, Christian Kreibich, Kirill Levchenko, Brandon Enright,
  Geoffrey~M. Voelker, Vern Paxson, and Stefan Savage, \emph{Spamalytics: {A}n
  empirical analysis of spam marketing conversion}, Communications of the ACM
  \textbf{52} (2009), no.~9, 99--107.

\bibitem[KSG13]{kosinski2013}
Michal Kosinski, David Stillwell, and Thore Graepel, \emph{Private traits and
  attributes are predictable from digital records of human behavior},
  Proceedings of the National Academy of Sciences of the United States of
  America \textbf{110} (2013), no.~15, 5802--5805.

\bibitem[LAH{\etalchar{+}}22]{lovato2022}
Juniper~L. Lovato, Antoine Allard, Randall Harp, Jeremiah Onaolapo, and Laurent
  H\'{e}bert-Dufresne, \emph{Limits of individual consent and models of
  distributed consent in online social networks}, Proceedings of the 2022 ACM
  Conference on Fairness, Accountability, and Transparency (New York Vity, NY,
  USA), FAccT '22, Association for Computing Machinery, 2022, pp.~2251--2262.

\bibitem[Les23]{leswing2023}
Kif Leswing, \emph{Microsoft limits {Bing A.I.} chats after the chatbot had
  some unsettling conversations}, CNBC (2023), 17 February, available at
  \url{https://www.cnbc.com/2023/02/17/microsoft-limits-bing-ai-chats-after-the-chatbot-had-some-unsettling-conversations.html}.

\bibitem[LHF{\etalchar{+}}21]{lazer2021}
David Lazer, Eszter Hargittai, Deen Freelon, Sandra Gonz\'alez-Bail\'on, Kevin
  Munger, Katherine Ognyanova, and Jason Radford, \emph{Meaningful measures of
  human society in the twenty-first century}, Nature \textbf{595} (2021),
  189--196.

\bibitem[LKG{\etalchar{+}}08]{ttt2008}
Kevin Lewis, Jason Kaufman, Marco Gonzalez, Andreas Wimmer, and Nicholas~A.
  Christakis, \emph{Tastes, ties, and time: {A} new (cultural, multiplex, and
  longitudinal) social network dataset using {F}acebook.com}, Social Networks
  \textbf{30} (2008), no.~4, 330--342.

\bibitem[LMP18]{loukides2018}
Mike Loukides, Hilary Mason, and {DJ} Patil, \emph{{Ethics and Data Science}},
  O'Reilly Media, Inc., Sebastopol, CA, USA, 2018.

\bibitem[Luc14]{luca2014}
Michael Luca, \emph{Were {OkCupid's} and {F}acebook's experiments unethical?},
  Harvard Business Review (2014), July, available at
  \url{https://hbr.org/2014/07/were-okcupids-and-facebooks-experiments-unethical}.

\bibitem[M\"{u}22]{muller2022}
Dennis M\"{u}ller, \emph{Situating {``Ethics in Mathematics''} as a philosophy
  of mathematics ethics education}, arXiv:2202.00705 (2022).

\bibitem[MC23]{muller2023b}
Dennis M\"{u}ller and Maurice Chiodo, \emph{Mathematical artifacts have
  politics: {T}he journey from examples to embedded ethics}, arXiv:2308.04871
  (2023).

\bibitem[McA16]{compton}
Nathan McAlone, \emph{{Why {`Straight Outta Compton'} had different {F}acebook
  trailers for people of different races}}, Business Insider (2016), available
  at
  \url{https://www.businessinsider.com/why-straight-outta-compton-had-different-trailers-for-people-of-different-races}.

\bibitem[MCF22]{hippo2021}
Dennis M\"{u}ller, Maurice Chiodo, and James Franklin, \emph{A {Hippocratic
  Oath} for mathematicians? {M}apping the landscape of ethics in mathematics},
  Science and Engineering Ethics \textbf{28} (2022), 41.

\bibitem[MCM09]{markines2009}
Benjamin Markines, Ciro Cattuto, and Filippo Menczer, \emph{Social spam
  detection}, Proceedings of the 5th International Workshop on Adversarial
  Information Retrieval on the Web (New York, NY, USA), AIRWeb '09, Association
  for Computing Machinery, 2009, pp.~41--48.

\bibitem[md4]{md4sg}
\emph{{Mechanism Design for Social Good}}, available at
  \url{https://www.md4sg.com} (accessed 31 December 2023).

\bibitem[MI]{MI-ox-data}
University of~Oxford Mathematical~Institute, \emph{Research using data
  involving humans}, available at
  \url{https://www.maths.ox.ac.uk/members/policies/data-protection/research-using-data-involving-humans}
  (accessed 26 December 2023).

\bibitem[MMS{\etalchar{+}}21]{mehrabi2021}
Ninareh Mehrabi, Fred Morstatter, Nripsuta Saxena, Kristina Lerman, and Aram
  Galstyan, \emph{A survey on bias and fairness in machine learning}, ACM
  Computing Surveys \textbf{54} (2021), no.~6, 115.

\bibitem[Mos21]{moses2021}
Bob Moses, \emph{Returning to `normal' in education is not good enough}, The
  Imprint (2021), 24 August, available at
  \url{https://imprintnews.org/opinion/returning-to-normal-in-education-is-not-good-enough/58069}.

\bibitem[MWW{\etalchar{+}}23]{moynihan2023}
Benjamin Moynihan, Robin~T. Wilson, Joan Wynne, Lee~J. McEwan, Mary~M. West,
  Frank~E. Davis, Herbert Clemens, Greg Budzban, Edith~Aurora Graf, and Aidan
  Soguero, \emph{{What’s math got to do with it?: Bob Moses, algebra, and the
  movement for civil rights, January 23, 1935--July 25, 2021}}, Notices of the
  American Mathematical Society \textbf{70} (2023), 258--277.

\bibitem[MWZ{\etalchar{+}}19]{mitchell2019-model-cards}
Margaret Mitchell, Simone Wu, Andrew Zaldivar, Parker Barnes, Lucy Vasserman,
  Ben Hutchinson, Elena Spitzer, Inioluwa~Deborah Raji, and Timnit Gebru,
  \emph{Model cards for model reporting}, Proceedings of the Conference on
  Fairness, Accountability, and Transparency (New York City, NY, USA), FAT*
  '19, Association for Computing Machinery, 2019, pp.~220--229.

\bibitem[N{\etalchar{+}}20]{ntoutsi2020}
Eirini Ntoutsi et~al., \emph{Bias in data-driven artificial intelligence
  systems---{An} introductory survey}, WIREs Data Mining and Knowledge
  Discovery \textbf{10} (2020), e1356.

\bibitem[Nar10]{arvind2010}
Arvind Narayanan, \emph{Cookies, supercookies and ubercookies: {S}tealing the
  identity of {W}eb visitors}, 33 Bits of Entropy: The End of Anonymous Data
  and What to Do About It (2010), available at
  \url{https://33bits.wordpress.com/2010/02/18/cookies-supercookies-and-ubercookies-stealing-the-identity-of-web-visitors/}.

\bibitem[{Nat}18]{proactive2018}
{National Academies of Sciences, Engineering, and Medicine}, \emph{Proactive
  policing: Effects on crime and communities}, The National Academies Press,
  Washington, DC, USA, 2018, available at \url{https://doi.org/10.17226/24928}.

\bibitem[Nob18]{noble2018}
Safiya~Umoja Noble, \emph{{Algorithms of Oppression: How Search Engines
  Reinforce Racism}}, NYU Press, New York City, NY, USA, 2018.

\bibitem[NS08]{arvind2008}
Arvind Narayanan and Vitaly Shmatikov, \emph{Robust de-anonymization of large
  sparse datasets}, 2008 IEEE Symposium on Security and Privacy (sp 2008),
  2008, pp.~111--125.

\bibitem[OBS{\etalchar{+}}19]{osoba2019}
Osonde~A. Osoba, Benjamin Boudreaux, Jessica Saunders, J.~Luke Irwin, Pam~A.
  Mueller, and Samantha Cherney, \emph{{Algorithmic Equity: A Framework for
  Social Applications}}, {RAND Corporation}, Santa Monica, CA, USA, 2019,
  available at \url{https://www.rand.org/pubs/research_reports/RR2708.html}.

\bibitem[{Off}22]{OSTP2022}
{Office of Science and Technology Policy, The White House}, \emph{{Blueprint
  for an AI Bill of Rights: Making Automated Systems Work for the American
  People}}, 2022, available at
  \url{https://www.whitehouse.gov/ostp/ai-bill-of-rights/}.

\bibitem[{OHR}]{ohrpp}
{OHRPP}, \emph{{UCLA Office of the Human Research Protection Program}},
  available at \url{https://ohrpp.research.ucla.edu} (accessed 26 December
  2023).

\bibitem[O'N16]{cathy2016}
Cathy O'Neil, \emph{{Weapons of Math Destruction: How Big Data Increases
  Inequality and Threatens Democracy}}, The Crown Publishing Group, New York,
  NY, USA, 2016.

\bibitem[PBL24]{purificato2024}
Erasmo Purificato, Ludovico Boratto, and Ernesto William~De Luca, \emph{User
  modeling and user profiling: {A} comprehensive survey}, arXiv:2402.09660
  (2024).

\bibitem[Pea24]{rat2024}
Jordan Pearson, \emph{Scientific journal publishes {AI}-generated rat with
  gigantic penis in worrying incident}, VICE (2024), 15 February, available at
  \url{https://www.vice.com/en/article/dy3jbz/scientific-journal-frontiers-publishes-ai-generated-rat-with-gigantic-penis-in-worrying-incident}.

\bibitem[Per11]{ttt-news}
Marc Perry, \emph{Harvard researchers accused of breaching students' privacy},
  The Chronicle of Higher Education (2011), available at
  \url{https://www.chronicle.com/article/harvard-researchers-accused-of-breaching-students-privacy/}.

\bibitem[Por21]{mason-slides}
Mason~A. Porter, \emph{{Data Ethics for Mathematicians}}, 2021, available at
  \url{https://zerodivzero.com/short_course/aaac8c66007a4d23a7aa14857a3b778c/title/181b207c7d2941278be4641ea5fe0e21}.

\bibitem[PSO{\etalchar{+}}20]{pierson2020}
Emma Pierson, Camelia Simoiu, Jan Overgoor, Sam Corbett-Davies, Daniel Jenson,
  Amy Shoemaker, Vignesh Ramachandran, Phoebe Barghouty, Cheryl Phillips, Ravi
  Shroff, and Sharad Goel, \emph{A large-scale analysis of racial disparities
  in police stops across the {United States}}, Nature Human Behaviour
  \textbf{4} (2020), 736--745.

\bibitem[RBKL20]{raghavan2020}
Manish Raghavan, Solon Barocas, Jon Kleinberg, and Karen Levy, \emph{Mitigating
  bias in algorithmic hiring: {E}valuating claims and practices}, Proceedings
  of the 2020 Conference on Fairness, Accountability, and Transparency (FAccT
  '20) (New York, NY, USA), Association for Computing Machinery, 2020,
  pp.~469--481.

\bibitem[Ree18]{royal2018}
Chris Reed, \emph{How should we regulate artificial intelligence?},
  Philosophical Transactions of the Royal Society A \textbf{376} (2018),
  no.~2128, 20170360.

\bibitem[RH23]{ross2023}
Casey Ross and Bob Herman, \emph{{Denied by AI: How Medicare Advantage} plans
  use algorithms to cut off care for seniors in need}, STAT (2023), available
  at
  \url{https://www.statnews.com/2023/03/13/medicare-advantage-plans-denial-artificial-intelligence/}.

\bibitem[RSMCM24]{rycroft2022}
Lucy Rycroft-Smith, Dennis M{\"u}ller, Maurice Chiodo, and Darren Macey,
  \emph{A useful ethics framework for mathematics teachers}, pp.~359--394,
  Springer International Publishing, Cham, Switzerland, 2024.

\bibitem[RT24]{AI-math}
Talia Ringer and Terrence Tao, \emph{{AI} for {M}ath resources}, available at
  \url{https://docs.google.com/document/d/1kD7H4E28656ua8jOGZ934nbH2HcBLyxcRgFDduH5iQ0/edit}
  (accessed 21 April 2024), 2024.

\bibitem[SAH{\etalchar{+}}21]{samba2021}
Nithya Sambasivan, Erin Arnesen, Ben Hutchinson, Tulsee Doshi, and Vinodkumar
  Prabhakaran, \emph{Re-imagining algorithmic fairness in {I}ndia and beyond},
  Proceedings of the 2021 ACM Conference on Fairness, Accountability, and
  Transparency (New York City, NY, USA), FAccT '21, Association for Computing
  Machinery, 2021, pp.~315--328.

\bibitem[Sal16]{salganik2016}
Matthew Salganik, \emph{{Computational Social Science: Social Research in the
  Digital Age}}, 2016, Sociology 596, Princeton University. Available at
  \url{https://www.princeton.edu/~mjs3/soc596_f2016/}.

\bibitem[Sal17]{salganik2017}
\bysame, \emph{{Bit by Bit: Social Research in the Digital Age}}, illustrated
  edition ed., Princeton University Press, Princeton, NJ, USA, 2017, available
  at \url{https://www.bitbybitbook.com}.

\bibitem[SBF{\etalchar{+}}19]{selbst2019}
Andrew~D. Selbst, Danah Boyd, Sorelle~A. Friedler, Suresh Venkatasubramanian,
  and Janet Vertesi, \emph{Fairness and abstraction in sociotechnical systems},
  Proceedings of the Conference on Fairness, Accountability, and Transparency
  (New York City, NY, USA), FAT* '19, Association for Computing Machinery,
  2019, pp.~59--68.

\bibitem[Sch21]{referee2021}
G.~William Schwert, \emph{The remarkable growth in financial economics,
  1974--2020}, Journal of Financial Economics \textbf{140} (2021), no.~3,
  1008--1046.

\bibitem[SH22]{stark2023}
Luke Stark and Jevan Hutson, \emph{Physiognomic artificial intelligence},
  Fordham Intellectual Property, Media and Entertainment Law Journal
  \textbf{32} (2022), no.~4, 2.

\bibitem[Sku21]{skufca2021}
Joe Skufca, \emph{Incorporating ethical discussions in the mathematics
  classroom}, SIAM News (2021), available at
  \url{https://sinews.siam.org/Details-Page/incorporating-ethical-discussions-in-the-mathematics-classroom-2}.

\bibitem[SM23]{stobaugh2023}
Ben Stobaugh and Dhiraj Murthy, \emph{Predicting gender and political
  affiliation using mobile payment data}, arXiv:2302.08026 (2023).

\bibitem[SMBM23]{straub2022}
Vincent~J. Straub, Deborah Morgan, Jonathan Bright, and Helen Margetts,
  \emph{Artificial intelligence in government: {C}oncepts, standards, and a
  unified framework}, Government Information Quarterly \textbf{40} (2023),
  no.~4, 101881.

\bibitem[SSGN17]{su2017}
Jessica Su, Ansh Shukla, Sharad Goel, and Arvind Narayanan,
  \emph{De-anonymizing {W}eb browsing data with social networks}, Proceedings
  of the 26th International Conference on World Wide Web (Republic and Canton
  of Geneva, CHE), WWW '17, International World Wide Web Conferences Steering
  Committee, 2017, pp.~1261--1269.

\bibitem[Ste22]{stemwedel2022}
Janet~D. Stemwedel, \emph{Science must not be used to foster white supremacy},
  Scientific American (2022), 24 May, available at
  \url{https://www.scientificamerican.com/article/science-must-not-be-used-to-foster-white-supremacy/}.

\bibitem[SVW21]{sadowski2021}
Jathan Sadowski, Salom\'{e} Viljoen, and Meredith Whittaker, \emph{Everyone
  should decide how their digital data are used --- {N}ot just tech companies},
  Nature \textbf{595} (2021), 169--171.

\bibitem[SWM05]{stumpf2005}
Michael P.~H. Stumpf, Carsten Wiuf, and Robert~M. May, \emph{Subnets of
  scale-free networks are not scale-free: {S}ampling properties of networks},
  Proceedings of the National Academy of Sciences of the United States of
  America \textbf{102} (2005), no.~12, 4221--4224.

\bibitem[Tao24]{tao-AI}
Terrence Tao, \emph{Two announcements: {AI for Math} resources, and
  \url{erdosproblems.com}}, available at
  \url{https://terrytao.wordpress.com/2024/04/19/two-announcements-ai-for-math-resources-and-erdosproblems-com/}
  (accessed 21 April 2024), 2024.

\bibitem[Tho]{thomas-tweet2020}
Rachel Thomas, available at
  \url{https://twitter.com/math_rachel/status/1297255819965169664} (accessed 26
  December 2023).

\bibitem[Tho21]{thomas-videos}
\bysame, \emph{11 short videos about {AI} ethics}, fast.ai (2021), available at
  \url{https://www.fast.ai/2021/08/16/eleven-videos/}.

\bibitem[TMP12]{Traud2012}
A.~L. Traud, P.~J. Mucha, and M.~A. Porter, \emph{Social structure of
  {F}acebook networks}, Physica A \textbf{391} (2012), no.~16, 4165--4180.

\bibitem[TPB24a]{tractenberg2024a}
Rochelle~E. Tractenberg, Victor Piercey, and Catherine~A. Buell, \emph{Defining
  ``ethical mathematical practice" through engagement with discipline-adjacent
  practice standards and the mathematical community}, Science and Engineering
  Ethics \textbf{30} (2024), 15.

\bibitem[TPB24b]{tractenberg2024b}
\bysame, \emph{{Proto Ethical Guidelines for Mathematical Practice}},
  https://osf.io/x5ur9/ (2024).

\bibitem[Uga17]{jugander2017medium}
Johan Ugander, \emph{Truth, lies, and an ethics of personalization}, Medium
  (2017), available at
  \url{https://medium.com/@jugander/truth-lies-and-an-ethics-of-personalization-e4ccfa7f2b84#.rzap3hm70}.

\bibitem[Uga20]{ugander2020}
\bysame, \emph{{Data Privacy and Data Ethics}}, 2020, Management Science and
  Engineering 234, Stanford University. Available at
  \url{https://web.stanford.edu/group/msande234/cgi-bin/wordpress/}.

\bibitem[Uga23]{Ugander2023}
\bysame, \emph{{Social Algorithms}}, 2023, Management Science and Engineering
  231, Stanford University. Available at
  \url{https://msande231.github.io/syllabus}.

\bibitem[{{Uni}}a]{datatheory}
{{University of California, Los Angeles}}, \emph{{Data Theory at UCLA}},
  available at \url{https://datatheory.ucla.edu} (accessed 26 December 2023).

\bibitem[{Uni}b]{ucla-AI}
{University of California, Los Angeles}, \emph{{UCLA: Generative AI}},
  available at \url{https://genai.ucla.edu/} (accessed 21 April 2024).

\bibitem[V\'{e}21]{veliz2021}
Carissa V\'{e}liz (ed.), \emph{{The Oxford Handbook of Digital Ethics}}, Oxford
  University Press, Oxford, UK, 2021.

\bibitem[Vin16]{hatebot2016}
James Vincent, \emph{Twitter taught {Microsoft's AI} chatbot to be a racist
  asshole in less than a day}, The Verge (2016), 24 March, available at
  \url{https://www.theverge.com/2016/3/24/11297050/tay-microsoft-chatbot-racist}.

\bibitem[VO23]{verma2023}
Pranshu Verma and Will Oremus, \emph{{ChatGPT} invented a sexual harassment
  scandal and named a real law prof as the accused}, The Washington Post
  (2023), 5 April, available at
  \url{https://www.washingtonpost.com/technology/2023/04/05/chatgpt-lies/}.

\bibitem[Wei21]{smbc2021}
Zach Weinersmith, \emph{Software}, Saturday Morning Breakfast Cereal (2021),
  available at \url{https://www.smbc-comics.com/comic/software}.

\bibitem[Whi12]{white2012}
Martha~C. White, \emph{Orbitz shows higher prices to {M}ac users}, Time (2012),
  26 June, available at
  \url{https://business.time.com/2012/06/26/orbitz-shows-higher-prices-to-mac-users/}.

\bibitem[Wika]{creative}
Wikipedia, \emph{{Creative Commons} license}, available at
  \url{https://en.wikipedia.org/wiki/Netflix_Prize#Cancelled_sequel} (accessed
  26 December 2023).

\bibitem[Wikb]{sequel}
\bysame, \emph{{Netflix Prize: Cancelled} sequel}, available at
  \url{https://en.wikipedia.org/wiki/Netflix_Prize#Cancelled_sequel} (accessed
  26 December 2023).

\bibitem[Wikc]{social-credit}
\bysame, \emph{{Social Credit System}}, available at
  \url{https://en.wikipedia.org/wiki/Social_Credit_System} (accessed 26
  December 2023).

\bibitem[Wikd]{gebru-google}
\bysame, \emph{{Timnit Gebru: Exit from Google}}, available at
  \url{https://en.wikipedia.org/wiki/Timnit_Gebru#Exit_from_Google} (accessed
  26 December 2023).

\bibitem[WJ23]{wiggins2023}
Chris Wiggins and Matthew~L. Jones, \emph{{How Data Happened: A History from
  the Age of Reason to the Age of Algorithms}}, W. W. Norton \& Company, New
  York City, NY, USA, 2023.

\bibitem[WSAO{\etalchar{+}}21]{wagner2021}
Claudia Wagner, Markus Strohmaier, Emre~Kıcıman Alexandra~Olteanu, Noshir
  Contractor, and Tina Eliassi-Rad, \emph{Measuring algorithmically infused
  societies}, Nature \textbf{595} (2021), 197--204.

\bibitem[WT24]{walk2024}
Stephen~M. Walk and Rochelle~E. Tractenberg, \emph{Helping students deal with
  ethical reasoning: {T}he proto-{G}uidelines for {E}thical {P}ractice in
  {M}athematics as a deck of cards}, arXiv:2403.16849 (2024).

\bibitem[WW24]{wiggins2024}
Madi Whitmann and Chris Wiggins, \emph{{Data: Past, Present, and Future}},
  2024, History-APMA UN2901, Columbia University. Available at
  \url{https://data-ppf.github.io} (accessed 31 December 2023).

\end{thebibliography}


\newcommand{\etalchar}[1]{$^{#1}$}
\providecommand{\bysame}{\leavevmode\hbox to3em{\hrulefill}\thinspace}
\providecommand{\MR}{\relax\ifhmode\unskip\space\fi MR }
\providecommand{\MRhref}[2]{%
  \href{http://www.ams.org/mathscinet-getitem?mr=#1}{#2}
}
\providecommand{\href}[2]{#2}


\end{document}